\documentclass[12pt,oneside]{article}
\newcommand{\Path}{}
\usepackage{ \Path Paper_style_2}
\usepackage{ \Path Keywords_WB}
\usepackage{ \Path Environments}
\usepackage{ \Path Category}

\usepackage{tikz-cd}
\usepackage{ \Path tikzcd2}
\usepackage{sidecap}
\usepackage{yfonts} 
\usepackage[most]{tcolorbox}

\DeclareMathOperator{\LN}{\text{\swabfamily LN}}
\DeclareMathOperator{\network}{\text{\swabfamily{Net}}}
\newcommand{\Time}{\mathcal{T}}
\DeclareMathOperator{\conn}{conn}
\DeclareMathOperator{\Bridge}{Bridge}
\begin{document}
	\title{Network of limit sets in dynamical systems}
	\author{Suddhasattwa Das\footnotemark[1]}
	\footnotetext[1]{Department of Mathematics and Statistics, Texas Tech University, Texas, USA}
	\date{\today}
	\maketitle
	\begin{abstract} 
	Points within a topological dynamical systems can be separated into transient and recurrent points. Most of the asymptotic properties of the dynamics is borne by the closure of the collection of recurrent points. This article presents a representation scheme for dynamics using the connectivity information in terms of the recurrent set. This representation is functorial and represents certain compositional features of dynamical systems.
	\end{abstract}
	
	\begin{keywords} Recurrence transient, omega-limit, functor, category \end{keywords}
	\begin{AMS}	37M99, 18F60, 37M22, 7B02, 37B35 \end{AMS}
	
	
	\section{Introduction} \label{sec:intro}
	
	While all dynamics is the result of a transition rule, such as a map or vector field, Dynamical systems theory itself is primarily concerned with the asymptotic properties of the dynamics \cite{Robinson1998dyn}. Asymptotic properties range from simple considerations such as periodic orbits to more complex behavior such as chaotic sets and invariant measures. A dynamical system should be viewed as a collection of numerous co-existing asymptotic behavior. This view point has been stressed and investigated in the context of reconstruction of dynamical systems \cite[e.g.]{BerryDas2023learning, BerryDas2024review, DasSuda2024recon}. One needs a way of representing a dynamical system which is not just its transition rule but also retains some information about some of these asymptotic properties.
	
	This article presents a representation of a dynamical system as a network in which the nodes are the various $\omega$ and $\alpha$ limit sets, and the connections are various subsets of the transient points in the dynamics. Dynamical systems have a \textit{categorical} structure \cite[e.g.]{DasSuda2025enrich, DasSuda2026join} from the various semi-conjugacy relations that bind different dynamics to each other. The representation we present is also shown to be \textit{functorial}, meaning that it translates these semi-conjugacies into a network of relations that exist within the structures we propose.
	
	There has been several other notable strategies for representing or encoding asymptotic properties. One of the earliest was Smale's \cite{PalisSmale1970structural} spectral decomposition of Axiom-A systems. A large class of these techniques are based on Conley's topological characterization of isolated invariant sets \cite{Conley1978} and the Morse decomposition of the dynamics which separates it into regions of gradient like dynamics, and non-wandering sets \cite{Conley1988_gradient, Conley1978}. These principles have been put to use in a variety of ways. The two most prominent applications is the theory of combinatorial vectors fields \cite[e.g.]{mrozek1990leray, Mrozek1996multivalued} as a means of representing the global properties of a flow, and the use of cohomology to determine the connections \cite[e.g.]{FranzosaMischaikow1988conn, BushEtAl2012comb, BanKalies2006comput, Cohen1995floer}. This extracts the algebraic information about how orbits leave the vicinity of an invariant set and enter the vicinity of another. Conley index theory has subsequently been extended to discrete time systems as well \cite[e.g.]{Kurland1981, RobbinSalamon1988shape}.
	
	A major drawback of these techniques is that it is unable to provide an insight into the internal structure of transitive dynamics. Transitivity is a property incompatible with gradient like flows. A gradient structure does not permit orbits to return near their initial points. However, this is one of the defining features of recurrence in ergodic theory, and the general chaotic systems. A chaotic system such as an non-uniformly hyperbolic system $X$ is densely packed with unstable periodic orbit. Not only the periodic points, but their stable manifolds, unstable manifolds, heteroclinic intersections are all dense in $X$. In a Morse decompositions, such a region $X$ would be represented by a single indivisible node. 
	
	Another means of representing a dynamical system by their asymptotic properties was the use of \emph{lattice-theory} by Kalies and others \cite[e.g.]{KaliesEtAl2014lattice, KaliesEtAl2016lattice, KaliesEtAl2021lattice}. Although this representation is well suited for computation \cite{KaliesEtAl2005algo} it is not well suited for the study of a single transitive invariant set. The absence of any proper sub-attractor leads to a trivial lattice without any interesting structure. 
	
	One of the features of our representation that distinguishes it from the others discussed so far is that although it relies on a graph-like structure, it is the edge set, not the node-set, which provides a partitioning of the dynamics. It is based on the simple principle that any point in a compact invariant set has non-empty $\omega$- and $\alpha$-limit sets. Every point thus acts as a connection between these two sets. The points in the state-space can be partitioned based on which pair of sets they connect. We call the elements of these partitions \emph{bridges}, and their limit sets their \textit{stations}.
	
	The biggest challenge in the utilization of this representation is in a topological realization of this network. This is partly helped by the fact that the network stores the poset structure of the limit-sets created by inclusion. The properties of stability and recurrence of an invariant set $A$ is a relative term, and depends upon which larger invariant set it is being considered to be a part of. While any point $x$ converges to $\omega(x)$ and every in point in $\omega(x)$ is a recurrent point in the orbit of $x$, the same would be not be true within $\omega(x)$ itself. Thus $\omega(x)$ might not be a transitive set, and the dynamics confined to $\omega(x)$ might reveal another layer of decomposition into limit sets, with points flowing from one to another. 
	
	A third strategy for retaining the asymptotic behavior of dynamics was through an outer approximation of the targeted invariant set \cite{Das2024zero, Das2024pipe}. The strategy here is to first approximate the deterministic dynamics as a low-noise Markov process. One can obtain an approximation of the topology of invariant set to any desired degree of accuracy, at the cost of losing determinism. The second stage is a topological realization of this Markov process as a deterministic skew-product flow, leading to an overall approximation that preserves both the dynamics and the transitivity. However, it loses sigh of the finer details within the dynamics, such as the nesting of invariant sets described above.
	
	An interesting outcome of the lattice structures present in the above two perspectives is the creation of a lattice of Lyapunov functions  \cite{RobbinSalamon1992lyap}. This lattice is collectively realized through a \emph{Lyapunov map}, a map from the underlying set of the dynamics into a simplicial complex created by the poset structure. This is also possible in the framework presented in this article.
	
	\section{Category}
	
	\paragraph{Categories} The categorical perspective provides a unified structural framework that reveals hidden connections and common patterns across seemingly unrelated fields. By abstracting away specific details to focus on relationships and transformations, it allows researchers to transfer powerful insights and solutions from one domain to another.
	
	\begin{definition}[Category]
		A category $\calC$ consists of a collection of objects $Ob(\calC)$ and, for every pair of objects $X, Y \in Ob(\calC)$, a collection of morphisms $\Hom_{\calC}(X, Y)$. The category is equipped with a composition operator $\circ: \Hom_{\calC}(Y, Z) \times \Hom_{\calC}(X, Y) \to \Hom_{\calC}(X, Z)$ defined for all morphisms $f: X \to Y$ and $g: Y \to Z$, as well as an identity morphism $\Id_X \in \Hom_{\calC}(X, X)$ for every object $X \in Ob(\calC)$. These must strictly satisfy the associativity axiom, $h \circ (g \circ f) = (h \circ g) \circ f$, and the identity laws, $f \circ \Id_X = f = \Id_Y \circ f$.
	\end{definition}
	
	\begin{definition}[Functor]
		A (covariant) functor $F: \calC \to \calD$ between categories $\calC$ and $\calD$ consists of a mapping of objects $X \mapsto F(X)$ and a mapping of morphisms $f \mapsto F(f)$, such that for every $f \in \Hom_{\calC}(X, Y)$, the image $F(f) \in \Hom_{\calD}(F(X), F(Y))$. Furthermore, $F$ must preserve the categorical structure by satisfying identity preservation, $F(\Id_X) = \Id_{F(X)}$, and composition preservation, $F(g \circ f) = F(g) \circ F(f)$.
	\end{definition}
	
	\begin{definition}[Natural Transformation]
		Let $F, G: \calC \to \calD$ be two functors. A natural transformation $\eta: F \Rightarrow G$ is a family of morphisms in $\calD$ indexed by the objects of $\calC$, such that for every $X \in Ob(\calC)$, there is a morphism $\eta_X: F(X) \to G(X)$ in $\calD$ (the components of $\eta$). This family must satisfy the commutativity condition: for every morphism $f: X \to Y$ in $\calC$, the diagram commutes, meaning $\eta_Y \circ F(f) = G(f) \circ \eta_X$.
	\end{definition}
	
	\begin{definition}[Dynamics as a Functor]
		Let $\Time$ be a monoid representing time (e.g., $(\mathbb{N}, +)$ for discrete time or $(\mathbb{R}_{\geq 0}, +)$ for continuous time). We view $\Time$ as a single-object category, $Ob(\Time) = \{\ast\}$, where $\Hom_{\Time}(\ast, \ast) = \Time$ and composition is the monoid operation. Let $\calC$ be a context category (e.g., $\Topo$, $\MeasCat$, or $\BanachCat$). A dynamical system over $\Time$ in the context $\calC$ is defined rigorously as a functor $\Phi: \Time \to \calC$. The object $X = \Phi(\ast) \in Ob(\calC)$ is the state space, and the morphism mapping $t \mapsto \Phi(t) \in \Hom_{\calC}(X, X)$ defines the time-evolution operators. The functor axioms automatically enforce the dynamical semigroup property: $\Phi(0) = \Id_X$ and $\Phi(t_1 + t_2) = \Phi(t_1) \circ \Phi(t_2)$.
	\end{definition}
	
	By a slight abuse of notation, the image $\Phi(t)$ will be denoted as $\Phi^t$. Formulating dynamics as a functor $\Phi: \Time \to \calC$ enables the application of universal properties, adjunctions, and categorical limits to system analysis and data-driven reconstruction. It formally decouples the underlying dynamics, the observation process, and the algorithmic mappings by treating them as distinct morphisms and natural transformations within a functor category $\Functor{\Time}{\calC}$. The most useful consequence is that the morphisms in this functor category, which are natural transformations, coincide exactly with \emph{semi-conjugacy} relations between dynamics :
	
	\begin{definition} [Semi-conjugacy]
		A semi-conjugacy between two dynamical systems $\paran{\Phi, X}$ and $\paran{\Phi', X'}$ is a continuous map $h:X\to X'$ satisfying the commutation :
		\begin{equation} \label{eqn:def:semiconjug}
			\begin{tikzcd}
				X \arrow[d, "h"'] \arrow[r, "\Phi^t"] & X \arrow[d, "h"] \\
				X' \arrow[r, "\Phi'^t"'] & X' \\
			\end{tikzcd} , \quad \forall t\in \Time.
		\end{equation}
		This is equivalent to a natural transformation $h:\Phi \Rightarrow \Phi'$ between the functors representing the two dynamics.
	\end{definition}
	
	\paragraph{Compact Hausdorff spaces} The category $\CompHaus$ of compact Hausdorff spaces possesses a uniquely rigid, algebraic structure that distinguishes it from general topological categories. Because its objects simultaneously support robust analytical tools (like regular Borel measures) and strict categorical properties, it serves as the natural domain for formulating continuous dynamical systems, probability theory, and operator algebras.
	
	\begin{lemma} [Compact Hausdorff spaces]
		\begin{enumerate}[(i)]
			\item $\CompHaus$ is isomorphic to the category of Eilenberg-Moore algebras for the ultrafilter monad on $\mathbf{Set}$.
			\item $\CompHaus$ has all small limits and colimits. Limits (such as products and pullbacks) coincide exactly with those constructed in $\Topo$. By Tychonoff's theorem, the arbitrary product of compact spaces remains compact.
			\item There is a categorical equivalence between the opposite category of compact Hausdorff spaces and the category of commutative, unital $C^*$-algebras, $\CompHaus \simeq \mathbf{cC}^*\mathbf{Alg}$.
			\item By the Riesz-Markov-Kakutani representation theorem, positive linear functionals on $C(X)$ correspond bijectively to regular Borel measures on $X$. This allows for a rigorous categorical treatment of Markov kernels and transition probabilities as morphisms within the Giry monad restricted to $\CompHaus$.
		\end{enumerate}
	\end{lemma}
	
	The first claim implies that despite being defined via topology, compact Hausdorff spaces are fundamentally algebraic structures. A compact Hausdorff space is perfectly determined by how ultrafilters of its underlying set converge, making the forgetful functor $U: \CompHaus \to \mathbf{Set}$ strictly monadic. The third claim implies a contravariant functor mapping $X \in \operatorname{Ob}(\CompHaus)$ to $C(X)$, the algebra of continuous complex-valued functions on $X$. This serves as the foundational bridge between topological dynamics and functional analysis, directly enabling Koopman operator frameworks where state-space evolution is dualized into operator evolution.
	
	\paragraph{A network category} We now build our category $\LN$. An object in this category is a poset $P$ with a terminal object $1_P$, along with a collection $\conn(P)$ of triples
	\[ \SetDef{ \paran{a,b,c} \in ob(\calP) \times ob(\calP) \times ob(\calP) }{ a\leq c,\, b\leq c }  .\]
	The objects of $\calP$ will be called \emph{stations} and $\calC$ will be called the \emph{connections}. The collection $\conn(P)$ must be such that every station has at least one outgoing or incoming bridge, Thus the objects can also be interpreted as directed multi-graphs whose underlying undirected simple graph is connected.  A morphism between two objects $P$ and $P'$ is a pair of mappings
	\[ \left.\begin{array}{c}
		f: \text{stations}(P) \to \text{stations}(P')  \\
		\tilde{f} : \text{bridges}(P) \to \text{bridges}(P') 
	\end{array} \right\} \,:\,
	\forall \paran{\begin{array}{c} a \\ b \\ c \end{array}} \in \text{bridges}(P) \,:\,
	\tilde{f} \paran{\begin{array}{c} a \\ b \\ c \end{array}} = \paran{\begin{array}{c} f(a) \\ f(b) \\ f(c) \end{array}}\]

	The connection of this combinatorial construct with dynamical systems is that there is a correspondence 
	\begin{tcolorbox}[colback=ChhaiF, colupper=black, colframe=black, coltext=kali, sharp corners, valign=center, top=0mm, bottom=0mm, left=0mm, right=0mm, boxrule=1.5pt]
		\begin{equation} \label{eqn:def:network}
			\begin{tikzcd} [scale cd = 2.0]
				\begin{array}{c} \Functor{\Time}{\CompHaus} \\ \mbox{Topological} \\ \mbox{dynamics} \end{array}
				\arrow[rr, Shobuj, "\network"] && 
				\begin{array}{c} \LN \\ \mbox{Bridges \&} \\ \mbox{stations} \end{array} 
			\end{tikzcd}
		\end{equation}
	\end{tcolorbox}
	which assigns every dynamical system $(\Phi, X)$ on a compact Hausdorff space $X$ its \emph{network of limit sets}. This is an object in $\LN$ defined as follows
	\begin{enumerate} [(i)]
		\item The poset $P$ is the collection of all $\alpha$- and $\omega$- limit sets, ordered by inclusion relation.
		\item Every point $x\in X$ has non-empty $\alpha(x)$ and $\omega(x)$. By reuse of notation, let $\alpha(x)$ and $\omega(x)$ also denote the corresponding nodes in $P$. If $x$ itself belongs to a limit set, $c(x)$ is defined to be the minimal node in $P$ whose underlying limit set contains $x$. Otherwise $c(x)$ is set to $1$, the terminal node in $P$. Thus one can associate to every $x\in X$ a triple
		\[\Bridge(x) := \paran{\alpha(x), \beta(x), c(x)} , \]
		of nodes from $P$.  The set $\conn(P)$ is defined to be all the collection of all such triples, disregarding multiplicity.
	\end{enumerate}
	
	The resulting network $\network(\Phi, X)$ thus represents both the poset structure of limit sets present in $X$, as well as the connectivity information. Connections for invariant sets are expresses using the language of stable and unstable sets. Recall that for any invariant subset $A$ :
	\[\begin{split}
		W^s(A) &:= \SetDef{x\in A}{ \lim_{t\to\infty} \dist \paran{ \Phi^t(x) , A } = 0 } , \\
		W^u(A) &:= \SetDef{x\in A}{ \lim_{t\to-\infty} \dist \paran{ \Phi^t(x) , A } = 0 } ,
	\end{split}\]
	For any node $c$ let $|c|$ denote the underlying subset of $X$. Thus any triple $(a,b,c)$ as shown on the left below :
	\[ (a,b,c) \in \conn \paran{\network(\Phi, X)} \imply \exists x\in |c| \,:\, x\in W^s(b) \cap W^u(a) . \]
	indicates that the stable and unstable sets of $b$ and $a$ respectively intersect through the orbit of a point in $|c|$. An elementary observation is that
	
	\begin{lemma}
		The points in $X$ are partitioned into the fibres of the mapping $\Bridge$. Thus the bridges of $\network(\Phi, X)$ represent a partitioning of the state-space $X$.
	\end{lemma}
	
	The most important realization is :
	
	\begin{theorem} [Functoriality of limit networks] \label{thm:LN_functor}
		The correspondence $\network$ from \eqref{eqn:def:network} is functorial. 
	\end{theorem}
	
	To verify, suppose one has a semi-conjugacy as shown in \eqref{eqn:def:semiconjug}. 
	
	\begin{lemma} \label{lem:LN_functor:1}
		Consider a semi-conjugacy as shown in \eqref{eqn:def:semiconjug} between topological dynamics. Then 
		\[ \forall x\in X, \;:\; h\paran{ \omega_{\Phi} (x) } = \omega_{\Phi'} \paran{ h(x) } . \]
	\end{lemma}
	
	Lemma \ref{lem:LN_functor:1} is proved in Section \ref{sec:proof:LN_functor:1}. The result can be restated as
	
	\begin{lemma} \label{lem:LN_functor:2}
		Consider a semi-conjugacy as shown in \eqref{eqn:def:semiconjug} between topological dynamics. Then if $A\subseteq X$ is invariant under $\Phi$ and $x\in A$, then $h(A)$ is invariant under $\Phi'$ and
		\[ \omega_{\Phi'} ( h(x) ) = h\paran{ \omega_{\Phi} (x) } , \quad \alpha_{\Phi'} ( h(x) ) = h\paran{ \alpha_{\Phi} (x) } . \]
	\end{lemma}
	
	A direct consequence of Lemma \ref{lem:LN_functor:2} is :
	
	\begin{lemma} \label{lem:LN_functor:3}
		Consider a semi-conjugacy as shown in \eqref{eqn:def:semiconjug} between topological dynamics. Then  for any $x,y\in X$ :
		\[ \Bridge_{\Phi}(x) = \Bridge_{\Phi}(y) \imply \Bridge_{\Phi'}(h(x)) = \Bridge_{\Phi'}(h(y)) . \]
	\end{lemma}
	
	Lemma \ref{lem:LN_functor:3} thus guarantees that each bridge in $\network(\Phi, X)$ is mapped into a unique bridge in $\network(\Phi', X')$. This leads to a creation of the functor $\LN$ as claimed.

	\begin{tcolorbox}[colback=ChhaiF, colupper=black, colframe=black, coltext=kali, sharp corners, valign=center, top=0mm, bottom=0mm, left=0mm, right=0mm, boxrule=1.5pt]
		\begin{equation}
			\begin{tikzcd} [scale cd = 2.0]
				\begin{array}{c} \Functor{\Time}{\CompHaus} \\ \mbox{Topological} \\ \mbox{dynamics} \end{array}
				\arrow[rr, Shobuj, "\network"] \arrow{dd}[name=n1, swap]{-} \arrow[dashed]{ddrr}[name=n2]{} && 
				\begin{array}{c} \LN \\ \mbox{Bridges \&} \\ \mbox{stations} \end{array} 
				\arrow[Itranga]{dd}{\text{Realize}} \\\\
				\Functor{\Time}{\CompHaus} && \Functor{\Time}{\CompHaus}
				\arrow[shorten <=1pt, shorten >=1pt, Rightarrow, Shobuj, to path={(n2) to[out=225,in=-45] node[yshift=15pt]{$\alpha$} (n1)} ]{  }
			\end{tikzcd}
		\end{equation}
	\end{tcolorbox}
	%
	
	\section{Recovering topological properties}

	Any dynamical system $\Phi^t : X \to X$ may be represented by such an object. Let $\calP$ be the collection of all subsets of $X$ which occur as $\omega$- or $\alpha$- limit sets. It inherits the partial ordering of the power set of $X$. Suppose there is a point $x$ such that $\omega(x) = B$ and $\alpha(x) = A$. Let $C$ be the minimal element of $\calP$ whose underlying set $\abs{C}$ contains $x$. Note that $C$ is an invariant set. Thus we have an orbit within the invariant set $C$ which connects $A$ to $B$. This information is stored as a triple $A, B, C$.
	
	\begin{definition} [Self-bridge]
		A bridge in an $\LN$ object which has the same source and destination stations, is called a \emph{self-bridge}.
	\end{definition}
	
	\begin{definition} [Recurrence]
		An $\LN$-object is called \emph{recurrent} if its terminal station has a self-bridge.
	\end{definition}
	
	\begin{lemma} [Transitivity of the dynamics]
		A dynamical system $(\Phi, X)$ is transitive iff 
	\end{lemma}
	
	\begin{lemma} [Transitivity of a subset]
		An asymptotic set $A$ of a dynamical system $(\Phi, X)$ is transitive iff there is an $\LN$-morphism $\mathbb{A} \to \network(\Phi, X)$ from a recurrent object $\mathbb{A}$, whose terminal station is mapped into the station corresponding to $A$.
	\end{lemma}
	
	\begin{definition} [Attractor]
		An invariant set $A$ of a dynamical system $(\Phi, X)$ is called an \textit{attractor} if its stable set contains a neighborhood of $A$.
	\end{definition}
	
	\begin{lemma} [Attractor / Repellor]
		Given a dynamical system $(\Phi, X)$, an invariant subset $A$ is an repellor if its represented by a station in $\network(\Phi, X)$ which has no incoming bridges. Similarly an attractor corresponds to a station in $\network(\Phi, X)$ which has no outgoing bridges. 
	\end{lemma}
	
	\begin{definition} [Attractor repellor pair]
		Given an attractor $A$ of a dynamical system $(\Phi, X)$, its \textit{dual repellor} is the set 
		\[ A^* := \SetDef{ x\in\Omega }{ \omega(x) \cap A = \emptyset } .\]
	\end{definition}
	
	\begin{lemma} [Attractor repellor pair]
		Let $A$ be a station in an $\LN$-object $\calN$. Suppose $A$ is an attractor, i.e. has no outgoing bridges. Then its dual repellor is the maximum of all those stations which have no bridge connecting to $A$.
	\end{lemma}
	
	\begin{definition} [Pseudo orbits]
		\(\{x_i\}\) (finite or infinite) is called a $\delta$-pseudo-orbit if the distance between the iterated point and the next point in the sequence satisfies:
		\(d(f(x_{i}),x_{i+1})<\delta \quad (\text{or\ }\varepsilon )\quad \text{for\ all\ }i\)
	\end{definition}
	
	\begin{definition} [Chain connectedness]
		Two points $x,y$ are said to be \emph{chain-connected} if for every $\delta>0$, there is a $\delta$-pseudo orbit starting at $x$ and ending at $y$.
	\end{definition}
	
	\begin{definition} [Chain transitivity graph]
		Being chain-connected is an transitive relation. The strongly connected components that emerge from this relation has the structure of a directed acyclic graph, called the \emph{chain transitivity graph}.
	\end{definition}
	
	\begin{definition} [Chain recurrent set]
		A strongly connected component in the chain transitivity graph is called a chain recurrent subset. 
	\end{definition}
	
	\begin{lemma} [Properties of chain recurrent sets]
		\begin{enumerate} [(i)]
			\item A chain recurrent component is a closed subset of $X$. 
			\item Restricting the flow or map to the chain recurrent set itself yields the same chain recurrent set.         
			\item Omega and alpha limit sets are contained within R(f).
			\item R(f) always contains the nonwandering set $\Omega(f)$. 
			\item The connected components of R(f) coincide with the chain-transitive components of the system, and each component is individually chain recurrent.
		\end{enumerate}
	\end{lemma}
	
	\section*{Declarations} 
	
	\paragraph{Funding} This work did not receive any funding or other monetary support.

\paragraph{Conflict of Interest Statement} – The author has no conflicts to disclose.

\paragraph{Ethics Approval Statement} There were no use of any living organisms in this research.

\paragraph{Author Contributions} The sole author SD was responsible for all of the research as well as in the preparation of the document.

\paragraph{Data Availability Statement} The research does not involve any experimental or proprietary data.

\paragraph{Declaration of generative AI use} All text, figures and research in this article was done manually by the author.
	\section{Appendix} \label{sec:appendix}
	
	\subsection{Proof of Lemma \ref{lem:LN_functor:1}} \label{sec:proof:LN_functor:1}
	
	Note that 
	\[\begin{split}
		y\in \omega_{\Phi} (x) &\,\Leftrightarrow\; \lim_{n\to\infty} \Phi^{t_n} x = y \,\Rightarrow\; \lim_{n\to\infty} h \paran{\Phi^{t_n} x} = h(y) \\
		&\,\Leftrightarrow\; \lim_{n\to\infty} \Phi'^{t_n} h(x) = h(y)\,\Leftrightarrow\; h(y) \in \omega_{\Phi'} (h(x)) .
	\end{split}\]
	This means that $ h\paran{\omega_{\Phi} (x)} \subseteq \omega_{\Phi'} (h(x))$. We next prove the reverse inclusion. 
	\[\begin{split}
		z\in \omega_{\Phi'} (h(x)) &\,\Leftrightarrow\; \lim_{n\to\infty} \Phi'^{t_n} h(x) = z \,\Leftrightarrow\; \lim_{n\to\infty} h\paran{\Phi'^{t_n} x} = z
	\end{split}\]
	Taking an appropriate subsequence we get 
	\[ \lim_{k\to\infty} \Phi'^{t_{n_k}} x = y \,\Leftrightarrow\; \lim_{k\to\infty} h \paran{\Phi'^{t_{n_k}} x} = h(y) . \]
	This means that $z=h(y)$. Thus every point in $\omega_{\Phi'} (h(x))$ is an image under $h$ of $\omega_{\Phi} (x)$. Thus $ h\paran{\omega_{\Phi} (x)} \supseteq \omega_{\Phi'} (h(x))$, proving the other half of the equality. \qed
	

\begin{thebibliography}{10}

\bibitem{Robinson1998dyn}
C.~Robinson.
\newblock \href{http://dx.doi.org/10.1201/9781482227871}{{\em Dynamical
  systems: stability, symbolic dynamics, and chaos}}.
\newblock CRC press, 1998.

\bibitem{BerryDas2023learning}
T.~Berry and S.~Das.
\newblock \href{http://dx.doi.org/10.1137/22M1516865}{Learning theory for
  dynamical systems}.
\newblock {\em SIAM J. Appl. Dyn.}, 22:2082 -- 2122, 2023.

\bibitem{BerryDas2024review}
T.~Berry and S.~Das.
\newblock
  \href{https://epubs.siam.org/doi/abs/10.1137/24M1696974?journalCode=siread}{Limits
  of learning dynamical systems}.
\newblock {\em SIAM review}, 16, 2025.

\bibitem{DasSuda2024recon}
S.~Das and T.~Suda.
\newblock \href{https://epubs.siam.org/doi/10.1137/25M1733203}{Dynamics, data
  and reconstruction}.
\newblock {\em Siam J. Appl. Dyn.}, 25, 2026.

\bibitem{DasSuda2025enrich}
S.~Das and T.~Suda.
\newblock \href{https://arxiv.org/pdf/2509.05900}{Dynamical systems as enriched
  functors}.
\newblock {\em arxiv}, 2509, 2025.

\bibitem{DasSuda2026join}
S.~Das and T.~Suda.
\newblock \href{https://arxiv.org/abs/2609.10912}{Joinings in {M}arkov
  categories}.
\newblock {\em arxiv}, 2609, 2026.

\bibitem{PalisSmale1970structural}
J.~Palis and S.~Smale.
\newblock \href{https://bookstore.ams.org/pspum-14}{Structural stability
  theorems}.
\newblock \href{https://bookstore.ams.org/pspum-14}{In {\em Global analysis}},
  volume~14, pages 223--231, 1970.

\bibitem{Conley1978}
C.~Conley.
\newblock {\em Isolated invariant sets and the Morse index}.
\newblock 38. Amer. Math. Soc., 1978.

\bibitem{Conley1988_gradient}
C.~Conley.
\newblock \href{http://dx.doi.org/10.1017/S0143385700009305}{The gradient
  structure of a flow: I}.
\newblock {\em Erg. Th. Dynam. Sys.}, 8(Charles Conley Memorial Issue), 1988.

\bibitem{mrozek1990leray}
M.~Mrozek.
\newblock \href{http://dx.doi.org/10.1090/S0002-9947-1990-0968888-1}{Leray
  functor and cohomological {C}onley index for discrete dynamical systems}.
\newblock {\em Trans. Amer. Math. Soc.}, 318(1):149--178, 1990.

\bibitem{Mrozek1996multivalued}
M.~Mrozek.
\newblock \href{http://dx.doi.org/10.1016/0898-1221(96)00127-7}{Topological
  invariants, mulitvalued maps and computer assisted proofs in dynamics}.
\newblock {\em Computers Math. Appl.}, 32(4):83--104, 1996.

\bibitem{FranzosaMischaikow1988conn}
R.~Franzosa and K.~Mischaikow.
\newblock \href{http://dx.doi.org/10.1016/0022-0396(88)90028-9}{The connection
  matrix theory for semiflows on (not necessarily locally compact) metric
  spaces}.
\newblock {\em J. Diff. Eq.}, 71(2):270--287, 1988.

\bibitem{BushEtAl2012comb}
J.~Bush et~al.
\newblock \href{http://dx.doi.org/10.1063/1.4772996}{Combinatorial-topological
  framework for the analysis of global dynamics}.
\newblock {\em Chaos}, 22(4), Dec 2012.

\bibitem{BanKalies2006comput}
H.~Ban and W.~Kalies.
\newblock \href{http://dx.doi.org/10.1115/1.2338654}{A computational approach
  to {C}onley's decomposition theorem}.
\newblock {\em Journal of Computational and Nonlinear Dynamics}, 1(4):312--319,
  2006.

\bibitem{Cohen1995floer}
R.~Cohen, J.~Jones, and G.~Segal.
\newblock Floer’s infinite dimensional morse theory and homotopy theory.
\newblock In {\em The Floer memorial volume}, pages 297--325. Springer, 1995.

\bibitem{Kurland1981}
H.~Kurland.
\newblock \href{http://dx.doi.org/10.1016/0022-0396(81)90028-0}{The {M}orse
  index of an isolated invariant set is a connected simple system}.
\newblock {\em J. Diff. Eq.}, 42(2):234--259, 1981.

\bibitem{RobbinSalamon1988shape}
J.~Robbin and D.~Salamon.
\newblock Dynamical systems, shape theory and the {C}onley index.
\newblock {\em Erg. Th. Dyn. Sys.}, 8(Charles Conley Memorial Issue):375--393,
  1988.

\bibitem{KaliesEtAl2014lattice}
W.~Kalies, K.~Mischaikow, and R.~Vandervorst.
\newblock \href{http://dx.doi.org/10.3934/jcd.2014.1.307}{Lattice structures
  for attractors i}.
\newblock {\em Journal of Computational Dynamics}, 1:307--338, 2014.

\bibitem{KaliesEtAl2016lattice}
W.~Kalies, K.~Mischaikow, and R.~Vandervorst.
\newblock \href{http://dx.doi.org/10.1007/s10208-015-9272-x}{Lattice structures
  for attractors ii}.
\newblock {\em Foundations of Computational Mathematics}, 16(5):1151--1191,
  2016.

\bibitem{KaliesEtAl2021lattice}
W.~Kalies, K.~Mischaikow, and R.~Vandervorst.
\newblock \href{http://dx.doi.org/10.1007/s10884-021-10056-8}{Lattice
  structures for attractors iii}.
\newblock {\em Journal of Dynamics and Differential Equations}, 34:1729--1768,
  2021.

\bibitem{KaliesEtAl2005algo}
W.~Kalies, K.~Mischaikow, and R.~Vandervorst.
\newblock \href{http://dx.doi.org/10.1007/s10208-004-0160-z}{An algorithmic
  approach to chain recurrence}.
\newblock {\em Foundations of Computational Mathematics}, 5(4):409--449, 2005.

\bibitem{Das2024zero}
S.~Das.
\newblock \href{http://dx.doi.org/10.48550/arXiv.2407.16673}{Reconstructing
  dynamical systems as zero-noise limits}.
\newblock {\em Tohoku Math. J.}, 77, 2025.

\bibitem{Das2024pipe}
S.~Das.
\newblock \href{doi.org/10.48550/arXiv.2407.16673}{Discrete-time dynamics,
  step-skew products, and pipe-flows}, 2024.

\bibitem{RobbinSalamon1992lyap}
J.~Robbin and D.~Salamon.
\newblock \href{http://dx.doi.org/10.1017/S0143385700006647}{Lyapunov maps,
  simplicial complexes and the {S}tone functor}.
\newblock {\em Erg. Th. Dyn. Sys.}, 12(1):153--183, 1992.

\end{thebibliography}

\end{document}